\documentclass[leqno,draft]{article}

\begin{document}

\title{A few remarks concerning \\ complex-analytic metric spaces}

\author{Stephen Semmes \\
        Rice University}

\date{}

\maketitle

        Let $E$ be a closed set in ${\bf C}^n$, normally with empty
interior, and let us consider continuously-differentiable functions on
$E$ in the sense of Whitney, in which the differential of the function
is automatically included.  This means a continuous complex-valued
function $f(z)$ on $E$ together with a real-linear mapping $df_z :
{\bf C}^n \to {\bf C}$ defined and continuous for $z \in E$, with
properties like those of an ordinary continuously-differentiable
function, in terms of local behavior.  Whitney's extension theorem
states that there is a continuously-differentiable function on ${\bf
C}^n$ in the usual sense which agrees with $f$ on $E$ and whose
differential is equal to $df_z$ when $z \in E$.  There are analogous
results for stronger smoothness properties, including $C^{1, \alpha}$
conditions.

        If $E$ is contained in a smooth real submanifold of ${\bf
C}^n$, then only the restriction of $df_z$ to the tangent space of the
submanifold is determined by $f$ on $E$.  However, the rest of the
differential is still relevant for continuously-differentiable
extensions to ${\bf C}^n$.  On a fractal set, it may be that $df_z$ is
determined in some or all directions by $f$ on $E$, even though $E$ is
not at all like a manifold in those directions.  This happens already
for self-similar Cantor sets, for instance.

        Every real-linear mapping from ${\bf C}^n$ into ${\bf C}$ can
be expressed in a unique way as the sum of a complex-linear mapping
and a conjugate-linear mapping.  For the differential $df_z$ of a
continuously-differentiable function $f$, these are typically denoted
$\partial f_z$ and $\overline{\partial} f_z$.  The condition
$\overline{\partial} f_z = 0$ for every $z \in E$ is a version of
holomorphicity for continuously-differentiable functions on $E$, which
is equivalent to the requirement that $df_z : {\bf C}^n \to {\bf C}$
be complex-linear for each $z \in E$.

        Of course, the restriction to $E$ of a holomorphic function on
a neighborhood of $E$ has this property.  It is already sufficient to
have a holomorphic extension in some directions, in such a way that
$df_z$ is the limit of complex-linear mappings for each $z \in E$.  If
$E$ is a complex submanifold of ${\bf C}^n$, then this condition
implies that $f$ is holomorphic on $E$ in the usual sense.  If $E$ is
totally disconnected, then there are plenty of nonconstant locally
constant functions $f$ on $E$, for which one can take $df_z = 0$ for
each $z \in E$.  There are also connected snowflake sets $E$ with a
lot of nonconstant continuously-differentiable functions with $df_z =
0$ for every $z \in E$.

        If $E$ is contained in a totally-real smooth submanifold of
${\bf C}^n$, then we can choose $df_z$ to be complex-linear, because a
real-linear mapping on a totally-real linear subspace of ${\bf C}^n$
can be extended to a complex-linear mapping on ${\bf C}^n$.
Similarly, if $E$ is contained in any smooth real submanifold, and
$df_z$ is complex-linear on the complex part of the tangent space of
the submanifold at $z$, then we can choose $df_z$ to be
complex-linear, since a real-linear mapping on a real-linear subspace
of ${\bf C}^n$ can be extended to a complex-linear mapping on ${\bf
C}^n$ exactly when it is complex-linear on the complex-linear part of
its domain.

        By contrast, even on a Cantor set, the whole differential
$df_z$ may be uniquely determined by $f$ on $E$, so that
complex-linearity of $df_z$ is a significant restriction.

        A remarkable phenomenon in several complex variables is that
holomorphic functions on some domains in ${\bf C}^n$, $n \ge 2$, can
always be extended to holomorphic functions on larger domains.  Related
results deal with holomorphic extensions in some directions of functions
on submanifolds that satisfy tangential Cauchy--Riemann equations.

        Ordinary holomorphic functions enjoy strong local regularity
properties, and in particular the limit of a sequence of holomorphic
functions on an open set that converges uniformly on compact subsets
is holomorphic.  If a sequence of holomorphic functions on a bounded
open set are continuous up to the boundary and converge uniformly on
the boundary, then they also converge uniformly on the interior by the
maximum principle.  Stability of tangential Cauchy--Riemann equations
in the $C^0$ topology can be analyzed on submanifolds using
integration by parts.

        On a Cantor set, arbitrary continuous functions can be
approximated uniformly by locally constant functions.  On snowflake
sets, continuous functions can be approximated in the $C^0$ topology
by continuously-differentiable functions with vanishing differential.

        Suppose that $E$ is a set with finite perimeter, and let ${\bf
1}_E$ be its characteristic function on ${\bf C}^n$.  If $f$ is a
continuously-differentiable function with $\overline{\partial} f_z =
0$ for each $z \in E$, then $\overline{\partial} (f \, {\bf 1}_E) = f
\, \overline{\partial} {\bf 1}_E$ in the sense of distributions.  This
equation makes sense for continuous functions $f$, and is preserved by
uniform convergence on compact sets.  A more classical version of this
in the complex plane is described in Exercise 2 at the end of Chapter
20 of \cite{r2}.  Note that the interior of $E$ may be empty.

        Let us restrict our attention now to compact sets $E$ in the
complex plane.  If the Lebesgue measure of $E$ is equal to $0$, then
every continuous function on $E$ can be approximated uniformly by the
restriction to $E$ of a holomorphic function on a neighborhood of $E$.
It suffices to start with the restriction to $E$ of a
continuously-differentiable function on ${\bf C}$ with compact
support, since these are dense among continuous functions on $E$.
Such a function can be represented as the Cauchy integral of its
$\overline{\partial}$ derivative, which can be approximated by the
Cauchy integral of the $\overline{\partial}$ derivative on the
complement of a neighborhood of $E$.

        Here are a couple of ways in which the condition
$\overline{\partial} b_z = 0$ for each $z \in E$ can be significant.
The kernel $(b(z) - b(w)) / (z - w)$ for $z, w \in E$ corresponds to
the commutator of multiplication by $b$ and the Cauchy integral
operator, and $\overline{\partial} b = 0$ on $E$ implies better
regularity for this kernel.  It also leads to smaller corrections for
the product of an extension of $b$ and a holomorphic function on a
domain with boundary $E$, say, to become holomorphic, by solving a
$\overline{\partial}$ problem.  It is even better if $b$ has
additional smoothness, which can also be seen in terms of extensions
whose $\overline{\partial}$ derivative vanishes more quickly near $E$.

        For the broader story of complex-analytic metric spaces, the
special case of subsets of ${\bf C}^n$ is instructive in the way that
there can be a lot of structure whether or not there is something like
ellipticity around.

\end{document}